\documentclass[10pt,a4paper]{article}
\usepackage{bbm}
\usepackage{mathrsfs}
\usepackage{amsfonts}
\usepackage{}

\usepackage{tikz}
\usetikzlibrary{positioning}
\usetikzlibrary{arrows}

\usepackage[leqno]{amsmath}
\usepackage{graphicx}
\usepackage{latexsym}
\usepackage{amsmath,amsfonts,amssymb,amsthm,mathrsfs,euscript,makeidx,color}
\usepackage{enumerate}
\usepackage{color}
\usepackage[colorlinks,linkcolor=blue,anchorcolor=green,citecolor=red]{hyperref}
\usepackage{indentfirst}
\renewcommand{\baselinestretch}{1.0}
\def\sqr#1#2{{\vcenter{\vbox{\hrule height.#2pt
              \hbox{\vrule width.#2pt height#1pt \kern#1pt \vrule width.#2pt}
              \hrule height.#2pt}}}}
\def\signed #1{{\unskip\nobreak\hfil\penalty50
              \hskip2em\hbox{}\nobreak\hfil#1
              \parfillskip=0pt \finalhyphendemerits=0 \par}}
\def\endpf{\signed {$\sqr69$}}

\def\5n{\negthinspace \negthinspace \negthinspace \negthinspace \negthinspace }
\def\4n{\negthinspace \negthinspace \negthinspace \negthinspace }
\def\3n{\negthinspace \negthinspace \negthinspace }
\def\2n{\negthinspace \negthinspace }
\def\1n{\negthinspace }

\def\dbR{\mathbb{R}}

\def\sK{\mathscr{K}}

\def\BN{{\bf N}}

\def\ds{\displaystyle}

\def\ns{\noalign{\ss}}
\def\no{\noindent}

\def\ss{\smallskip}
\def\ms{\medskip}

\def\q{\quad}
\def\qq{\qquad}

\def\({\Big (}
\def\){\Big )}
\def\[{\Big[}
\def\]{\Big]}

\def\lan{\langle}
\def\ran{\rangle}

\def\rf{\eqref}

\def\b{\beta}
\def\g{\gamma}
\def\d{\delta}

\def\z{\zeta}
\def\k{\kappa}

\def\m{\mu}
\def\n{\nu}

\def\t{\tau}
\def\f{\varphi}
\def\th{\theta}
\def\o{\omega}

\def\i{\infty}

\def\h{\widehat}

\def\cd{\cdot}

\def\les{\leqslant}
\def\ges{\geqslant}
\def\bde{\begin{definition}\label}
\def\ede{\end{definition}}
\def\be{\begin{equation}}
\def\bel{\begin{equation}\label}
\def\ee{\end{equation}}
\def\bt{\begin{theorem}\label}
\def\et{\end{theorem}}
\def\bc{\begin{corollary}\label}
\def\ec{\end{corollary}}
\def\bl{\begin{lemma}\label}
\def\el{\end{lemma}}
\def\bp{\begin{proposition}\label}
\def\ep{\end{proposition}}
\def\bas{\begin{assumption}\label}
\def\eas{\end{assumption}}
\def\br{\begin{remark}\label}
\def\er{\end{remark}}
\def\bex{\begin{example}\label}
\def\ex{\end{example}}
\def\ba{\begin{array}}
\def\ea{\end{array}}
\def\ben{\begin{enumerate}}
\def\een{\end{enumerate}}

\def\square#1{\vbox{\hrule\hbox{\vrule height#1%
     \kern#1\vrule}\hrule}}
\def\rectangle#1#2{\vbox{\hrule\hbox{\vrule height#1%
     \kern#2\vrule}\hrule}}

\font\tenbb=msbm10 \font\sevenbb=msbm7 \font\fivebb=msbm5

\newfam\bbfam
\scriptscriptfont\bbfam=\fivebb \textfont\bbfam=\tenbb
\scriptfont\bbfam=\sevenbb

\newtheorem{theorem}{\indent Theorem}[section]
\newtheorem{definition}[theorem]{\indent Definition}
\newtheorem{proposition}[theorem]{\indent Proposition}
\newtheorem{corollary}[theorem]{\indent Corollary}
\newtheorem{lemma}[theorem]{\indent Lemma}
\newtheorem{remark}[theorem]{\indent Remark}
\newtheorem{example}[theorem]{\indent Example}

\newtheorem{assumption}[theorem]{\indent Assumption}

\makeatletter
   
   \@addtoreset{equation}{section}
\makeatother

\begin{document}

\title{\bf Remarks on Viscosity Super-Solutions of Quasi-Variational Inequalities}

\author{Yue Zhou\footnote{School of Mathematics and Statistics, Central South University, Changsha, Hunan, China. Email:{\tt zhouyuemath@csu.edu.cn.}},
~~~Xinwei Feng\footnote{Corresponding author. Zhongtai Securities Institute for Financial Studies, Shandong University, Jinan, Shandong
250100, China. This author was supported by National Natural Science
Foundation of China (No. 12001317),  Shandong Provincial Natural Science Foundation (No. ZR2020QA019) and QILU Young Scholars Program of Shandong University. Email:{\tt xwfeng} {\tt @sdu.edu.cn.}},~~~Jiongmin Yong\footnote{Department of
Mathematics, University of Central Florida, Orlando, FL 32816, USA. This author was supported in part by NSF Grant DMS-1812921. Email:{\tt jiongmin.yong@ucf.edu.} }}

\maketitle

\no\bf Abstract: \rm For Hamilton-Jacobi-Bellman (HJB) equations, with the standard definitions of viscosity super-solution and sub-solution, it is known that there is a comparison between any (viscosity) super-solutions and sub-solutions. This should be the same for HJB type quasi-variational inequalities (QVIs) arising from optimal impulse control problems. However, according to a natural adoption of the definition found in \cite{Barles 1985, Barles 1985b}, the uniqueness of the viscosity solution could be guaranteed, but the comparison between viscosity super- and sub-solutions could not be guaranteed. This short paper introduces a modification of the definition for the viscosity super-solution of HJB type QVIs so that the desired comparison theorem will hold.

\ms

\no\bf Keywords: \rm Hamilton-Jacobi-Bellman quasi-variational inequality, viscosity solution, comparison theorem, optimal impulse control.

\ms

\no\bf AMS Mathematics Subject Classification. \rm 49N25, 49L20, 49L25.

\section{Introduction}

For standard optimal control problem of ordinary differential equations in finite time horizon, Bellman's dynamic programming method leads to the following Hamilton-Jacobi-Bellman (HJB, for short) equation:
\bel{HJB1}\left\{\2n\ba{ll}
\ds V_t(t,x)+H(t,x,V_x(t,x))=0,\qq(t,x)\in[0,T]\times\dbR^n,\\
\ns\ds V(T,x)=h(x),\qq x\in\dbR^n.\ea\right.\ee
The value function $V(\cd\,,\cd)$ of the optimal control problem is characterized as the unique viscosity solution to the above HJB equation.
\ms

For optimal impulse control problem of ordinary differential equations, it is known that, under some proper conditions, instead of HJB equation \rf{HJB1}, the value function happens to be the unique viscosity solution to the following HJB type quasi-variational inequality (QVI, for short):
\bel{QVI1}\left\{\2n\ba{ll}
\ds\min\big\{V_t(t,x)+H(t,x,V_x(t,x)),\BN[V](t,x)-V(t,x)\big\}=0,\q(t,x)\in[0,T]\times\dbR^n,\\
\ns\ds V(T,x)=h(x),\qq x\in\dbR^n,\ea\right.\ee
where
$$\BN[V](t,x)=\inf_{\xi\in K}\(V(t,x+\xi)+\ell(t,\xi)\),\qq(t,x)\in[0,T]\times\dbR^n$$
is called an obstacle operator, for some given function $\ell(\cd\,,\cd)$, called the impulse cost, and $K\subseteq\dbR^n$ is some convex closed cone in which the impulse control takes values (see \cite{Bensoussan-Lions 1973, Bensoussan-Lions 1984}).

\ms

We now recall the following definitions, for convenience of presentation (\cite{Crandall-Lions 1983, Barron-Evans-Jensen 1984, Ishii 1984, Crandall-Ishii-Lions 1992}).

\bde{HJB-vs} \rm (i) A continuous function $V:[0,T]\times\dbR^n\to\dbR$ is called a
{\it viscosity sub-solution} of HJB equation \rf{HJB1} if the following holds
\bel{V<h}V(T,x)\les h(x),\qq\forall x\in\dbR^n,\ee
and for any smooth function $\f(\cd\,,\cd)$ as long as $V(\cd\,,\cd)-\f(\cd\,,\cd)$ attains a local maximum at $(t_0,x_0)\in[0,T)\times\dbR^n$, the following holds:
\bel{f_t+H>0}\f_t(t_0,x_0)+H(t_0,x_0,\f_x(t_0,x_0))\ges0.\ee

(ii) A continuous function $V:[0,T]\times\dbR^n\to\dbR$ is called a {\it viscosity super-solution} of HJB equation \rf{HJB1} if the following holds
\bel{V>h}V(T,x)\ges h(x),\qq\forall x\in\dbR^n,\ee
and for any smooth function $\f(\cd\,,\cd)$ as long as $V(\cd\,,\cd)-\f(\cd\,,\cd)$ attains a local minimum at $(t_0,x_0)\in[0,T)\times\dbR^n$, the following holds:
\bel{f_t+H<0}\f_t(t_0,x_0)+H(t_0,x_0,\f_x(t_0,x_0))\les0.\ee

(iii) In the case that $V(\cd\,,\cd)$ is both viscosity sub- and super-solution of HJB equation \rf{HJB1}, it is called a {\it viscosity solution} of HJB equation \rf{HJB1}.

\ede

The following is an adoption of those in \cite{Barles 1985, Barles 1985b} for QVI \rf{QVI1}, see also \cite{Barron-Evans-Jensen 1984,Tang-Yong 1993}.

\bde{QVI-vs} \rm (i) A continuous function $V:[0,T]\times\dbR^n\to\dbR$ is called a {\it viscosity sub-solution} of QVI \rf{QVI1} if \rf{V<h} holds,
and for any smooth function $\f(\cd\,,\cd)$ as long as $V(\cd\,,\cd)-\f(\cd\,,\cd)$ attains a local maximum at $(t_0,x_0)\in[0,T)\times\dbR^n$, the following holds:
\bel{f_t+H*>0}\min\Big\{\f_t(t_0,x_0)+H(t_0,x_0,\f_x(t_0,x_0)),\BN[V](t_0,x_0)
-V(t_0,x_0)\Big\}\ges0.\ee

(ii) A continuous function $V:[0,T]\times\dbR^n\to\dbR$ is called a {\it viscosity super-solution} of QVI \rf{QVI1} if \rf{V>h} holds,
and for any smooth function $\f(\cd\,,\cd)$ as long as $V(\cd\,,\cd)-\f(\cd\,,\cd)$ attains a local minimum at $(t_0,x_0)\in[0,T)\times\dbR^n$, the following holds:
\bel{f_t+H*<0}\min\Big\{\f_t(t_0,x_0)+H(t_0,x_0,\f_x(t_0,x_0)),\BN[V](t_0,x_0)
-V(t_0,x_0)\Big\}\les0.\ee

(iii) In the case that $V(\cd\,,\cd)$ is both viscosity sub- and super-solution of QVI \rf{QVI1}, it is called a viscosity solution of QVI \rf{QVI1}.

\ede

It is well-known that under some mild conditions, if $V(\cd\,,\cd)$ and $\h V(\cd\,,\cd)$ are respectively viscosity sub-solution and viscosity super-solution of HJB equation \rf{HJB1}, then the following comparison holds:
\bel{V<V}V(t,x)\les\h V(t,x),\qq(t,x)\in[0,T]\times\dbR^n.\ee
This implies the uniqueness of viscosity solution to HJB equation \rf{HJB1} (see, for example \cite{Crandall-Lions 1983, Ishii 1984, Crandall-Ishii-Lions 1992, Bardi-Capuzzo-Dolcetta 1997}). On the other hand, for QVI \rf{QVI1}, although it is known that (see \cite{Tang-Yong 1993}) under proper conditions, the viscosity solution is unique, it seems that a general comparison between viscosity sub-solution and viscosity super-solution is missing, to out best knowledge. It turns out that in order to establish such a general comparison theorem, we have to modify the above definition of viscosity super-solution. This is the main purpose of this short paper.

\section{Modification of Viscosity Super-Solution Definition for QVIs}

For convenience of presentation, let us first introduce the following assumptions.

\ms

{\bf(H1)} Let $h:\dbR^n\to\dbR$ and $H:[0,T]\times\dbR^n\times\dbR^n\to\dbR$ be continuous and the following hold:
\bel{H-H}\ba{ll}
\ns\ds h(x)\ges-h_0,\qq |h(x)-h(y)|\leq L(1+|x|^\mu+|y|^\mu)|x-y|^\gamma,\\
\ns\ds|H(t,x,p)|\les L(1+|x|^\m)(1+|p|),\\
\ns\ds|H(t,x,p)-H(s,y,q)|\les L(1+|x|^\m\vee|y|^\m)|p-q|+\o\big(|x|\vee|y|+|p|\vee|q|,|t-s|+|x-y|),\\
\ns\ds\qq\qq\qq\qq\qq\qq\qq\qq\qq\qq\qq\forall t\in[0,T],~x,y,p,q\in\dbR^n,\ea\ee
for some constants $h_0,L>0$, $\m,\gamma\in[0,1)$ and continuous function $\o:[0,\i)\times[0,\i)\to[0,\i)$, monotone increasing in each argument and $\o(r,0)=0$ for any $r\ges0$.

\ms

{\bf(H2)} Let $\ell:[0,T]\times K\to(0,\i)$ be continuous with $K\subset\dbR^n$ being a closed convex cone. There exist constants $\ell_0>0$ such that
\bel{ell(1)}\ell_0\les\ell(t,\xi),\qq(t,\xi)\in[0,T]\times K,\ee
\bel{ell(3)}\ba{ll}
\ns\ds\ell(t,\xi+\xi')\les\ell(t,\xi)+\ell(t,\xi')
,\qq t\in[0,T],~\xi,\xi'\in K,\ea\ee
\bel{ell(4)}\ell(t',\xi)\les\ell(t,\xi),\qq t\les t',~\xi\in K.\ee

\bel{ell(5)}\lim_{|\xi|\rightarrow\infty}\inf_{t\in[0,T]}\ell(t,\xi)=+\infty.\ee

With Definition \ref{QVI-vs} (i), we have the following result.

\bl{Lemma-5.1} \sl A continuous function $V(\cd\,,\cd)$ is a viscosity sub-solution of \rf{QVI1} if and only if the following holds
\bel{V<N[V]}V(t,x)\les\BN[V](t,x),\qq\forall(t,x)\in[0,T]\times\dbR^n,\ee
and it is a viscosity sub-solution of the following HJB equation:
\bel{HJB-V0}\left\{\2n\ba{ll}
\ds V_t(t,x)+H\big(t,x,V_x(t,x)\big)=0,\q
(t,x)\in[0,T)\times\dbR^n,\\
\ns\ds V(T,x)=h(x),\qq x\in\dbR^n.\ea\right.\ee

\el

\it Proof. \rm Let $V(\cd\,,\cd)$ be a viscosity sub-solution of \rf{QVI1}. If \rf{V<N[V]} is false, then there exists a point $(t_0,x_0)\in[0,T]\times\dbR^n$ such that
$$V(t_0,x_0)>\BN[V](t_0,x_0),$$
and by the continuity of the involved functions, there exists a $\d_0>0$ such that
\bel{v(t,x)}V(t,x)\ges\BN[V](t,x)+\d_0,\qq\forall |t-t_0|+|x-x_0|\les\d_0,~(t,x)\in[0,T]\times\dbR^n.\ee
Let $\z(\cd\,,\cd)\in C^\i([0,T],\dbR^n)$ such that
$$\left\{\2n\ba{ll}
\ns\ds\z(t_0,x_0)=1,\\
\ns\ds0\les\z(t,x)<1,\qq\forall(t,x)\in[0,T]\times\dbR^n,~(t,x)\ne(t_0,x_0),\\
\ns\ds\z(t,x)=0,\qq\qq\forall(t,x)\in[0,T]\times\dbR^n,~|t-t_0|+|x-x_0|\ges\d_0.\ea\right.$$
For $R>\ds\max_{(t,x)\in[0,T]\times\dbR^n,|t-t_0|+|x-x_0|\les\d_0}|V(t,x)|$, define
$$\phi(t,x)=V(t,x)+2R\zeta(t,x).$$
Then, for $(t,x)\in[0,T]\times\dbR^n$, $|t-t_0|+|x-x_0|=\d_0$, we have that
$$\phi(t,x)=V(t,x)<R\les V(t_0,x_0)+2R=\phi(t_0,x_0).$$
Therefore, there exists a point $(t_1,x_1)\in[0,T]\times\dbR^n$, $|t_1-t_0|+|x_1-x_0|<\d_0$ such that
$$\phi(t_1,x_1)\ges\max_{(t,x)\in[0,T]\times\dbR^n,|t-t_0|+|x-x_0|\les\d_0}\phi(t,x),$$
which means that $(t,x)\mapsto V(t,x)+2R\z(t,x)$ attains a local maximum at $(t_1,x_1)$.
By the definition of viscosity sub-solution, we have
$$\min\Big\{-2R\z_t(t_1,x_1)+H(t_1,x_1,-2R\z_x(t_1,x_1)),\BN[V](t_1,x_1)-V(t_1,x_1)
\Big\}\ges0.$$
Hence
$$V(t_1,x_1)\les\BN[V](t_1,x_1),$$
which contradicts \rf{v(t,x)}. Next, let $\f(\cd\,,\cd)$ be smooth such that $V(\cd\,,\cd)-\f(\cd\,,\cd)$ attains a local maximum at $(t_0,x_0)\in[0,T]\times\dbR^n$. Then, according to Definition \ref{QVI-vs} (i), one has \rf{f_t+H*>0} which implies
$$\f_t(t_0,x_0)+H\big(t_0,x_0,\f_x(t_0,x_0)\big)\ges0.$$
Thus, $V(\cd\,,\cd)$ is a viscosity sub-solution of \rf{HJB-V0}. Conversely, if $V(\cd\,,\cd)$ is a viscosity sub-solution of \rf{HJB-V0} and \rf{V<N[V]} holds, then combination of these two implies \rf{f_t+H*>0}. Hence, $V(\cd\,,\cd)$ is a viscosity sub-solution of \rf{QVI1}. \endpf

\ms

We point out that if $V(\cd\,,\cd)$ is just a viscosity super-solution of \rf{QVI1} in the sense of Definition \ref{QVI-vs} (ii), we could not get \rf{V<N[V]}, and we do not know if it is a viscosity super-solution of \rf{HJB-V0}, either. In fact, if $V(\cd\,,\cd)-\f(\cd\,,\cd)$ attains a local minimum at $(t_0,x_0)\in[0,T]\times\dbR^n$, then \rf{f_t+H*<0} holds,
%
%
which could not exclude the following:
\bel{V>V>V}\f_t(t_0,x_0)+H\big(t_0,x_0,\f_x(t_0,x_0)\big)>0>\BN[V](t_0,x_0)-V(t_0,x_0).\ee
Therefore, with such a definition, if $\h V(\cd\,,\cd)$ is a viscosity super-solution of \rf{QVI1} (in the sense of Definition \ref{QVI-vs} (ii)) and $V(\cd\,,\cd)$ is a viscosity sub-solution of \rf{QVI1}, it is not clear whether the comparison \rf{V<V} is true.
The main reason is that \rf{f_t+H*<0} could not guarantee \rf{V<N[V]} (at $(t_0,x_0)$).

\ms

Note that in the standard proof of the uniqueness of viscosity solution to \rf{QVI1} (see \cite{Tang-Yong 1993}, as well as \cite{Barles 1985, Barles 1985b}), a key is that one has two possible viscosity solutions of \rf{QVI1} to compare (not the case that one is a sub-solution
and the other is a super-solution). Because both are already viscosity sub-solutions, by Lemma \ref{Lemma-5.1} (i), both satisfy the constraint \rf{V<N[V]} which cannot be guaranteed for a viscosity super-solution $\h V(\cd\,,\cd)$ in the sense of Definition \ref{QVI-vs} (ii). This is the major reason that \rf{V<V} cannot be proved if $\h V(\cd\,,\cd)$ is just a viscosity super-solution in the sense of Definition \ref{QVI-vs} (ii).

\ms

We now take a closer look at \rf{QVI1}. To reveal the essential feature, let us first look at HJB equation \rf{HJB-V0}. Since the unknown function $V(\cd\,,\cd)$ is not necessarily differentiable, one uses the notion of the viscosity sub- and/or super-solution in the following manner: Find a smooth function $\f(\cd\,,\cd)$ touching $V(\cd\,,\cd)$ at some point $(t_0,x_0)$ from below (or above), then in the equation, replace the derivatives of $V(\cd\,,\cd)$ by those of $\f(\cd\,,\cd)$, with the equality replaced by an inequality of suitable direction. Now, for QVI \rf{QVI1}, one can try to use the same idea. However, \rf{V>V>V} tells us that the constraint \rf{V<N[V]} cannot be guaranteed from the viscosity super-solution defined by Definition \ref{QVI-vs} (ii), due to the nonlocal nature of $\BN[V](\cd\,,\cd)$. Realizing this, we may modify the definition as follows.

\bde{super-solution} \rm A function $V(\cd\,,\cd)\in C([0,T]\times\dbR^n)$ is called a {\it viscosity super-solution} of QVI \rf{QVI1} if \rf{V>h}, \rf{V<N[V]} hold, and, for any smooth function $\f(\cd\,,\cd)$, as long as $V(\cd\,,\cd)-\f(\cd\,,\cd)$ attains a local minimum at $(t_0,x_0)\in[0,T)\times\dbR^n$ with
\bel{V(t_0)}V(t_0,x_0)<\BN[V](t_0,x_0),\ee
it holds
\bel{f_t<0}\f_t(t_0,x_0)+H\big(t_0,x_0,\f_x(t_0,x_0)\big)\les0.\ee

\ede

From \rf{V(t_0)}--\rf{f_t<0}, we see that constraint \rf{V<N[V]} is assumed first of all, and on the set $\{(t,x)\in[0,T]\times\dbR^n\bigm|V(t,x)=\BN[V](t,x)\}$, condition like \rf{f_t<0} is not required for the viscosity super-solution of QVI \rf{QVI1}. Thus, unlike Lemma \ref {Lemma-5.1}, viscosity super-solution of HJB equation \rf{HJB-V0} with constraint \rf{V<N[V]} is a viscosity super-solution of \rf{QVI1} in the sense of Definition \ref{super-solution}. But the converse might not be true due to the fact that \rf{V>V>V} might happen. Namely, a
viscosity super-solution of QVI \rf{QVI1} in the sense of Definition \ref{super-solution} might not be a viscosity super-solution of HJB \rf{HJB-V0}. We formally state this fact as follows.

\bl{Lemma-6.3} \sl Let $V(\cd\,,\cd)$ be a viscosity super-solution of HJB equation \rf{HJB-V0} satisfying constraint \rf{V<N[V]}. Then it is a viscosity super-solution of \rf{QVI1} in the sense of Definition \ref{super-solution}.

\el

The above definition of viscosity super-solution is stronger than that  Definition \ref{QVI-vs} (ii), i.e., $V(\cd\,,\cd)$ is a viscosity super-solution of QVI \rf{QVI1} in the sense of Definition \ref{super-solution}, then it is a viscosity super-solution of QVI \rf{QVI1} in the sense of Definition \ref{QVI-vs} (ii). More importantly, with the above stronger definition, the standard results that the value function of the corresponding optimal impulse control being a viscosity solution of the corresponding QVI \rf{QVI1} remains(for readers' convenience, we will give the proof in the Appendix), and by strengthening the definition, the uniqueness of viscosity solution remains as well.  The following example illustrates that the two definitions of viscosity super-solution are not equivalent.

\bex{} \rm Consider the following QVIs:
\bel{HJB-ex}\left\{\2n\ba{ll}
\ds\min\Big\{V_t(t,x)-V_x(t,x)+g(t,x),\BN[V](t,x)-V(t,x)\Big\}=0,\qq(t,x)\in[0,T]\times\dbR,\\
\ns\ds V(T,x)=h(x)\equiv xe^{-x},\qq x\in\dbR,\ea\right.\ee
where $g:[0,T]\times\dbR\to[0,\i)$ is continuous, and
$$\BN[V](t,x)=\inf_{\xi\ges0}\big\{V(t,x+\xi)+\ell_0+\ell_0\xi\big\},\qq(t,x)\in[0,T]\times\dbR,$$
with $\ell_0>0$ small, undetermined. Let
$$V(t,x)=(x-T+t)e^{-(x-T+t)},\qq(t,x)\in[0,T]\times\dbR.$$
Then
$$V(T,x)=xe^{-x},\qq x\in\dbR,$$
and
\bel{sub}V_t(t,x)-V_x(t,x)+g(t,x)=g(t,x)\ges0,\qq(t,x)\in[0,T]\times\dbR.\ee
Next, for fixed $t_0\in[0,T]$, we let $x_0=T-t_0+1\in\dbR$. Then $x_0-T+t_0=1$.
Note that
\bel{xe<e}xe^{-x}\les e^{-1},\qq\forall x\in\dbR.\ee
Thus, for $(t_0,x_0)=(t_0,T-t_0+1)$, one has
$$V(t_0,x_0+\xi)+\ell(\xi)=(x_0+\xi-T+t_0)e^{-(x_0+\xi-T+t_0)}+\ell_0(1+\xi)
=(1+\xi)[e^{-(1+\xi)}+\ell_0]\equiv\psi(\xi).$$
To find a minimum of $\psi(\cd)$, we set
$$0=\psi'(\xi)=[e^{-(1+\xi)}+\ell_0]-(1+\xi)e^{-(1+\xi)}=\ell_0
-\xi e^{-(1+\xi)}=e^{-1}[\ell_0e-\xi e^{-\xi}].$$
Now, noting \rf{xe<e}, if we choose $\ell_0<e^{-2}$, then there are
$0<\xi_1<1<\xi_2$ such that
$$\psi'(\xi)\left\{\2n\ba{ll}
\ns\ds>0,\qq\xi\in[0,\xi_1)\cup(\xi_2,\i),\\
\ns\ds<0,\qq\xi\in(\xi_1,\xi_2),\\
\ns\ds=0,\qq\xi=\xi_1,\xi_2.\ea\right.$$
Hence,
$$\BN[V](t_0,x_0)-V(t_0,x_0)=(1+\xi_2)[e^{-(1+\xi_2)}+\ell_0]
-e^{-1}.$$
Since $\xi_2>1$, we must have
$$(1+\xi_2)e^{-(1+\xi_2)}<e^{-1}.$$
Consequently, by shrinking $\ell_0>0$ further if necessary, we will have
$$\BN[V](t_0,x_0)-V(t_0,x_0)<0.$$
In addition, by the continuity, we can find $\d>0$ such that
$$\BN[V](t,x)-V(t,x)<0,\qq|t-t_0|+|x-x_0|\les\d,~(t,x)\in[0,T]\times\dbR.$$
Now, we let
$$g(t,x)=0,\qq(t,x)\notin\Big\{(t,x)\in[0,T]\times\dbR\bigm||t-t_0|+|x-t_0|<{\d\over2}
\Big\}.$$
Then
$$\min\Big\{V_t(t,x)-V_x(t,x)+g(t,x),\BN[V](t,x)-V(t,x)\Big\}\les0,$$
in the classical sense. According to Definition \ref{QVI-vs} (ii),
$V(\cd\,\cd)$ is a viscosity super-solution of \rf{HJB-ex}. However, since
constraint \rf{V<N[V]} is not satisfied, it is not a viscosity super-solution of
\rf{HJB-ex} in the sense of Definition \ref{super-solution}. From \rf{sub}, we see that $V(\cd\,,\cd)$ is actually a sub-solution of the HJB equation
\bel{HJB-ex*}\left\{\2n\ba{ll}
\ds V_t(t,x)-V_x(t,x)+g(t,x)=0,\qq(t,x)\in[0,T]\times\dbR,\\
\ns\ds V(T,x)=h(x)\equiv xe^{-x},\qq x\in\dbR,\ea\right.\ee
Hence, being able to prove a comparison between a viscosity sub-solution and a viscosity super-solution of \rf{QVI1} defined by Definition \ref{QVI-vs} essentially amounts to being able to prove a comparison between two viscosity sub-solutions of \rf{HJB-ex*} which is unlikely possible.\\
\indent Let us introduce the following sets. For function $u:[0,t]\times\dbR^n\rightarrow\dbR$, $(t,x)\in[0,T)\times\dbR^n$,
\begin{equation*}\begin{aligned}
\mathcal P^{1,+}u(t,x)=\{&(b,p)\in\dbR\times\dbR^n\bigm|u(t,x)\les u(s,y)+b(t-s)+\langle p,x-y\rangle\\
&+o(|t-s|+|x-y|),\text{ as }(t,x)\rightarrow(x,y)\},
\end{aligned}\end{equation*}
\begin{equation*}\begin{aligned}
\bar{\mathcal P}^{1,+}u(t,x)=\{&(b,p)\in\dbR\times\dbR^n\bigm|\exists(s_i,y_i)\in[0,T]\times\dbR^n,
(b_i,p_i)\in\mathcal P^{1,+}u(s_i,y_i),\\
& (s_i,y_i,u(s_i,y_i),b_i,p_i)\rightarrow(s,y,u(s,y),b,p)\}.
\end{aligned}\end{equation*}
We define
$$\mathcal P^{1,-}u(t,x)=-\mathcal P^{1,+}(-u)(t,x),\qq\bar{\mathcal P}^{1,-}u(t,x)=-\bar{\mathcal P}^{1,+}(-u)(t,x).$$
{Using the above notations, we have the following equivalent definition of viscosity sub- and (modified) super-solution.
\bde{Viscosity solution-equivalent definition} \rm (i) A continuous function $V(\cd\,,\cd)$ is called a {\it viscosity sub-solution} of HJB quasi-variational inequality \rf{QVI1} if
\bel{VT)<-equi}V(T,x)\les h(x),\qq x\in\dbR^n,\ee
and for any $(b,p)\in\bar{\mathcal P}^{1,+}V(t,x)$, it holds
\bel{min>0-equi}\min\Big\{b+H\big(t,x,p\big),\BN[V](t,x)-V(t,x)\Big\}\ges0.\ee
(ii) A function $V(\cd\,,\cd)\in C([0,T]\times\dbR^n)$ is called a {\it viscosity super-solution} of \rf{QVI1} if
\bel{VT)>-equi}V(T,x)\ges h(x),\qq x\in\dbR^n,\ee
\rf{V<N[V]} hold, and, for any $(b,p)\in\bar{\mathcal P}^{1,-}V(t,x)$ with
\bel{V(t_0)-equi}V(t,x)<\BN[V](t,x),\ee
it holds
\bel{f_t<0-equi}b+H\big(t,x,p\big)\les0.\ee
\ede}
\ex

\section{A Comparison Theorem}

Having modified the definition of viscosity super-solution, we will now present a general comparison theorem for viscosity sub- and super-solutions of QVIs. More precisely, besides \rf{QVI1}, we consider
\bel{QVI2}\left\{\ba{ll}
\ds\min\Big\{\bar V_t(t,x)+\bar H(t,x,\bar V_x(t,x)),\bar\BN[\bar V](t,x)-\bar V(t,x)
\Big\}=0,\q (t,x)\in[0,T)\times\dbR^n,\\
\ns\ds\bar V(T,x)=\bar h(x),\qq x\in\dbR^n,\ea\right.\ee
where $\bar H:[0,T]\times\dbR^n\times\dbR^n\to\dbR$ and $\bar h:\dbR^n\to\dbR$ are continuous functions also satisfying (H1), and $\bar\BN$ is defined as before with
$\ell(\cd\,,\cd\,,\cd)$ replaced by $\bar\ell(\cd\,,\cd\,,\cd)$ for which (H2) holds. Our general comparison theorem can be stated as follows.

\bt{comparison} \sl Let {\rm(H1)--(H2)} hold for $(h,H,\ell)$ and $(\bar h,\bar H,\bar\ell)$. Let
\bel{comp}\left\{\2n\ba{ll}
\ds -h_0\les h(x)\les\bar h(x),\qq\forall x\in\dbR^n,\\
\ns\ds H(t,x,p)\les\bar H(t,x,p),\qq\forall(t,x,p)\in[0,T]\times\dbR^n\times\dbR^n,\\
\ns\ds\ell(t,\xi)\les\bar\ell(t,\xi),\qq\forall(t,\xi)\in[0,T]
\times K.\ea\right.\ee
Let $V$ be a viscosity sub-solution of QVI \rf{QVI1}, and $\bar V$ be a viscosity super-solution of QVI \rf{QVI2}, satisfying the following:
\bel{prop}\left\{\2n\ba{ll}
\ds|V(t,x)|,|\bar V(t,x)|\les C(1+|x|^\g),\qq\forall(t,x)\in[0,T]\times
\dbR^n,\\
\ns\ds|V(t,x)-V(s,x')|,|\bar V(t,x)-\bar V(s,x')|\les C(|t-s|+|x-x'|^\k),\q\forall
t,s\in[0,T],~x,x'\in\dbR^n,\ea\right.\ee
for some $\g\in[0,1)$ and $\k\in(0,\b)$. Then
\bel{V<hV}V(t,x)\les\bar V(t,x),\qq(t,x)\in[0,T]\times\dbR^n.\ee
\et

\it Proof. \rm
We prove this by contradiction. Suppose $V\leq\bar V$ is not true. Then there exists a point
 $(\bar t,\bar x)\in(0,T)\times\dbR^n$
such that
$$V(\bar t,\bar x)-\bar V(\bar t,\bar x)\equiv \eta>0.$$
We take constant $G>0$ large enough and $\theta,\rho,\lambda\in(0,1]$ small enough, so that the following hold
\begin{equation}
\left\{
             \begin{array}{lr}
            \theta G<1,\q  G\ell_0-2(\frac{2C}{\alpha})^{\frac{1}{\beta-\kappa}}>0,\q\rho T<\frac{\eta}{4},   \\
            \theta G(V(\bar t,\bar x)+h_0)<\frac{\eta}{4},\q 2\theta \langle \bar x \rangle+\frac{\lambda}{\bar t}<\frac{\eta}{4}. \notag
             \end{array}
\right.
\end{equation}
We define
$$\Phi(t,x,y)=(1-\theta G)V(t,x)-\bar V(t,y)-\theta \frac{2T-t}{2T}(\langle x\rangle+\langle y\rangle)-\frac{1}{2\varepsilon}|x-y|^2+\rho t-\frac{\lambda}{t},$$where
$\langle x\rangle:=\sqrt{|x|^2+1}$.
It is clear that there exists $(t_0,x_0,y_0)\in(0,T)\times\dbR^n\times\dbR^n$ such that
 \begin{align}\label{sup*}
\Phi(t_0,x_0,y_0)&=\max_{t,x,y}\Phi(t,x,y)\ges\Phi(T,0,0)\notag\\
&=(1-\theta G)V(T,0)-\bar V(T,0)+\rho T-\frac{\lambda}{T}.\notag
\end{align}
This implies that for some (absolute) constant $C>0$,
$$\theta (\langle x_0\rangle+\langle y_0\rangle)+\frac{1}{\varepsilon}|x_0-y_0|^2\les C(1+|x_0|^\gamma+|y_0|^\gamma).$$
Thus, there exists a constant $C_\theta$, depending on $\theta$ such that
$$|x_0|,|y_0|,\frac{1}{\varepsilon}|x_0-y_0|^2\les C_\theta.$$
On the other hand, from $2\Phi(t_0,x_0,y_0)\ges \Phi(t_0,x_0,x_0)+\Phi(t_0,y_0,y_0)$, we have
$$\frac{1}{2\varepsilon}|x_0-y_0|^2\les C|x_0-y_0|^\gamma \rightarrow 0 \q as \q\varepsilon \rightarrow 0.$$
Next we claim that for $\varepsilon>0$ small, we have $t_0<T$. In fact, if $t_0=T$, then
\begin{align}
\Phi(\bar t,\bar x,\bar x)&\equiv(1-\theta G)V(\bar t,\bar x)-\bar V(\bar t,\bar x)-\theta \frac{2T-\bar t}{T}\langle \bar x\rangle+\rho \bar t-\frac{\lambda}{\bar t}\notag\\
&\les \Phi(t_0,x_0,y_0)\notag\\
&=(1-\theta G)h(x_0)-\bar h(y_0)-\frac{\theta}{2}(\langle x_0\rangle+\langle y_0\rangle)-\frac{1}{2\varepsilon}|x_0-y_0|^2+\rho T-\frac{\lambda}{T}\notag\\
&\les L(1+|x_0|^\mu+|y_0|^\mu)|x_0-y_0|^\gamma+\rho T+\theta Gh_0.\notag
\end{align}
Then, we have
\begin{align}
\eta&=V(\bar t,\bar x)-\bar V(\bar t,\bar x)\notag\\
&\les\theta GV(\bar t,\bar x)+L(1+|x_0|^\mu+|y_0|^\mu)|x_0-y_0|^\gamma+2\theta \langle \bar x \rangle+\rho T+\frac{\lambda}{\bar t}+\theta Gh_0\notag\\
&< \frac{3\eta}{4}+L(1+|x_0|^\mu+|y_0|^\mu)|x_0-y_0|^\gamma.\notag
\end{align}
Thus, for the given $\theta>0$,if $\varepsilon>0$ is small enough, the above gives a contradiction.\\
Next, we claim that
\begin{equation}\label{bar vleq}\bar V(t_0, y_0)<\bar\BN[\bar V](t_0, y_0),\end{equation}
Actually, if we have
$$\bar V(t_0, y_0)=\bar\BN[\bar V](t_0,y_0)=\bar V(t_0, y_0+\xi_0)+\bar\ell(t_0, y_0,\xi_0)$$
for some $\xi_0\in K$, then we have
$$\alpha|\xi_0|^\beta\les \bar\ell(t_0, y_0,\xi_0)\les C(1+|\xi_0|^\k))\Rightarrow|\xi_0|\les (\frac{2C}{\alpha})^{\frac{1}{\beta-\k}}.$$
Thus, as $\varepsilon\rightarrow 0$ , we have $\o\(|x_0|\vee|y_0|,|x_0-y_0|\)\rightarrow 0$, then
$$\ba{ll}
\ns\ds \Phi(t_0, x_0+\xi_0, y_0+\xi_0)-\Phi(t_0,x_0,y_0)\\
\ns\ds\ges(1-\theta G)[V(t_0, x_0+\xi_0)-V(t_0, x_0)]-[\bar V(t_0, y_0+\xi_0)-\bar V(t_0, y_0)]-\theta \frac{2T-t_0}{2T}(\langle x_0+\xi_0\rangle-\langle x\rangle+\langle y_0+\xi_0\rangle-\langle y_0\rangle)\\
\ns\ds\ges\theta G\ell(t_0,x_0,\xi_0)+\ell(t_0,y_0,\xi_0)-\ell(t_0,x_0,\xi_0)-2\theta|\xi_0|\\
\ns\ds\ges \theta G\ell_0-\o\(|x_0|\vee|y_0|,|x_0-y_0|\)-2\theta(\frac{2C}{\alpha})^{\frac{1}{\beta-\kappa}}\\
\ns\ds\ges \theta[G\ell_0-2(\frac{2C}{\alpha})^{\frac{1}{\beta-\kappa}}]-\o\(|x_0|\vee|y_0|,|x_0-y_0|\)>0.\ea$$
This contradicts the definition of $(x_0,y_0)$, hence \eqref{bar vleq} hold.\\

Now, let us define\\
$$\varphi(t,x,y)=\theta (\frac{2T-t}{2T})(\langle x\rangle+\langle y\rangle)+\frac{1}{2\varepsilon}|x-y|^2-\rho t+\frac{\lambda}{t},$$
then
\begin{equation}
\left\{
             \begin{array}{lr}
            \varphi_t(t,x,y)=-\rho-\frac{\theta}{2T}(\langle x\rangle+\langle y\rangle)-\frac{\lambda}{t^2},   \\
             \varphi_x(t,x,y)=\theta (\frac{2T-t}{2T})\frac{x}{\langle x\rangle}+\frac{1}{\varepsilon}(x-y),  \\
             \varphi_y(t,x,y)=\theta (\frac{2T-t}{2T})\frac{y}{\langle y\rangle}+\frac{1}{\varepsilon}(y-x). \notag
             \end{array}
\right.
\end{equation}
{Applying Theorem 9 of \cite{M.G.-H 1990} to the function
$$(1-\theta G)V(t,x)+(-\bar V(t,y))-\varphi(t,x,y)$$
at $(t_0,x_0,y_0)$, we can find $b,c\in \mathbb R$, $X,Y\in S(n)$  such that
\begin{equation*}\left\{\begin{aligned}
&(b,\varphi_x(t_0,x_0,y_0))\in\bar{\mathcal P}^{1,+}((1-\theta G)V)(t_0,x_0),\\
&(c,\varphi_y(t_0,x_0,y_0))\in\bar{\mathcal P}^{1,+}(-\bar V)(t_0,y_0),\\
&b+c=\varphi_t(t_0,x_0,y_0).
\end{aligned}\right.\end{equation*}}
By the definition of viscosity sub/sup-solution(Definition \ref{Viscosity solution-equivalent definition}), we have
$$-b-(1-\theta G)H(t_0,x_0,\frac{1}{1-\theta G}\varphi_x(t_0,x_0,y_0))\les 0,$$
$$c-\bar H(t_0,y_0,-\varphi_y(t_0,x_0,y_0))\ges 0.$$
Thus we have
\begin{align}
\rho+\frac{\theta}{2T}(\langle x_0\rangle+\langle y_0\rangle)&\les (1-\theta G)H(t_0,x_0,\frac{1}{1-\theta G}\varphi_x(t_0,x_0,y_0))-\bar H(t_0,y_0,-\varphi_y(t_0,x_0,y_0))\notag\\
&\les (1-\theta G)H(t_0,x_0,\frac{1}{1-\theta G}\varphi_x(t_0,x_0,y_0))- H(t_0,y_0,-\varphi_y(t_0,x_0,y_0))\notag\\
&=[L(1+|x_0|^\m\vee|y_0|^\m)\Big|{1\over1-\th G}\f_x(t_0,x_0,y_0)+\f_y(t_0,
x_0,y_0)\Big|]\notag\\
&\q+\o\(|x_0|\vee|y_0|+\big\{{1\over1-\th G}
|\f_x(t_0,x_0,y_0)|\big\}\vee|\f_y(t_0,x_0,y_0)|,|x_0-y_0|\)\notag\\&\q+\th GL(1+|x_0|^\m){1\over1-\th G}\(1+|\f_x(t_0,x_0,y_0)|\)\notag.
\end{align}
letting $\varepsilon\rightarrow0$, we have $(t_0,x_0,y_0)\rightarrow(t_0,x_0,x_0)$, thus
$$\rho\les L(1+|x_0|^\m)\({1\over1-\th G}+1\){\th(2 T-t_0)\over2 T}+{\th G L\over1-\th G}(1+|x_0|^\m)\(1+{\th(2 T-t_0)\over2 T}\)-\frac{\theta}{T}\langle x_0\rangle.$$
Finally, since $\mu<1$, by sending $\theta\rightarrow0$, we get $\rho\les0$ which is a contradiction.\endpf

\section{Appendix}

In this appendix, we will briefly present an optimal impulse control problem that leads to the quasi-variational inequality of form \rf{QVI1}. Consider the following equation:
\bel{state1}X(s)=x+\int_t^sf(r,X(r))dr+\xi(s),\qq s\in[t,T],\ee
where $f:[0,T]\times\dbR^n\to\dbR^n$ is a given map, $(t,x)\in[0,T)\times\dbR^n$ is called an {\it initial pair}, and
\bel{xi}\xi(s)=\sum_{k\ges1}\xi_k{\bf1}_{[\t_k,T]}(s),\qq s\in[t,T]\ee
is called an {\it impulse control} with $\{\t_k\}_{k\ges1}\subset[t,T]$ being a non-decreasing finite sequence, and $\xi_k\in K$, $k\ges1$, called {\it admissible impulses}, for some non-empty closed convex cone $K\subseteq\dbR^n$ with the vertex at the origin. In the above, we allow $\t_k=\t_{k+1}$ for some $k\ges1$. Let $\sK[t,T]$ be the set of all impulse controls of form \rf{xi}. Under some mild conditions, for any initial pair $(t,x)\in[0,T)\times\dbR^n$ and impulse control $\xi(\cd)\in\sK[t,T]$, equation \rf{state1} admits a unique solution $X(\cd)\equiv X(\cd\,;t,x,\xi(\cd))$. Clearly, both $\xi(\cd)$ and $X(\cd)$ are right-continuous. To measure the performance of the impulse control $\xi(\cd)$, we introduce the following cost functional
\bel{cost1}J(t,x;\xi(\cd))=\int_t^Tg(s,X(s))ds+h(X(T))+\sum_{k\ges1}\ell\big(\t_k,\xi_k\big),\ee
where
\bel{ghl}g:[0,T]\times\dbR^n\to[0,\infty),\q h:\dbR^n\to[0,\infty),\q\ell:[0,T]\times K\to(0,\infty)\ee
are suitable maps. Here, the terms on the right-hand side of \rf{cost1} are called the {\it running cost}, the {\it terminal cost} and the {\it impulse cost}, respectively. In the above, we may assume that $g$ and $h$ are just bounded uniformly from below. By a possible translation, we can simply assume that they are non-negative, for convenience. We emphasize that the impulse cost $\ell(t,\xi)$ is strictly positive. Mimicking the classical case, we formulate the following optimal impulse control problem.

\ms

{\bf Problem (C).} \rm For any initial pair $(t,x)\in[0,T)\times\dbR^n$, find a $\bar\xi(\cd)\in\sK[t,T]$ such that
\bel{infJ}J(t,x;\bar\xi(\cd))=\inf_{\xi(\cd)\in\sK[t,T]}J(t,x;\xi(\cd))=V(t,x).\ee

\ms

We call $\bar\xi(\cd)$ an {\it optimal impulse control}, the corresponding
$\bar X(\cd)\equiv X(\cd\,;t,x,\bar\xi(\cd))$ an {\it optimal state trajectory}, $(\bar X(\cd),\bar\xi(\cd))$ an {\it optimal pair}, and $V(\cd\,,\cd)$ the {\it value function} of Problem (C).

\ms

Before going further, let us first introduce the following hypotheses.

\ms



{\bf(H1)$'$} The map $f:[0,T]\times\dbR^n\to\dbR^n$ is continuous and there exists a constant $L>0$ such that
\bel{|f-f|}|f(t,x)-f(t,x')|\les L|x-x'|,\qq\forall t\in[0,T],~x,x'\in\dbR^n.\ee
\bel{|f|}|f(t,0)|\les L,\qq\forall t\in[0,T].\ee

\ms

{\bf(H1)$''$} Maps $g:[0,T]\times\dbR^n\to[0,\infty)$ and $h:\dbR^n\to[0,\infty)$ are continuous. There exist constants $L,\m>0$ and $0<\d\les1$ such that
\bel{g,h}0\les g(t,x),h(x)\les L\big(1+|x|^{\m+\d}\big),\qq\forall(t,x)\in[0,T]\times\dbR^n,\ee
\bel{g-g,h-h}|g(t,x)-g(t,x')|,|h(x)-h(x')|\les L\(1+|x|^\m\vee|x'|^\m\)|x-x'|^\d,\qq\forall t\in[0,T],x,x'\in\dbR^n.\ee

\ms

%
%
%
%


The above (H1)$'$--(H1)$''$ implies (H1). We will also let (H2) hold for the impulse cost $\ell(\cd\,,\cd)$. Condition \rf{ell(1)} implies that as long as an impulse is made, no matter how small the $\xi$ is, there is a strictly positive fixed cost $\ell_0$.  Condition \rf{ell(3)} means that if at $(t,x)$ an impulse of size $\xi+\xi'$ needs to be made, then one should make just one impulse of that size instead of making an impulse of size $\xi$ immediately followed by another with size $\xi'$. Hence, in an optimal impulse control, $\t_k<\t_{k+1}$ if both are impulsive moments.  Because of this condition, $\xi\mapsto\ell(t,\xi)$ should be ``sublinear''.  Condition \rf{ell(4)} means that if an impulse is going to be made, then the later the better, which is essentially due to the discount effect. In what follows, we call the impulse control that contains no impulses the {\it trivial impulse control}, denote it by $\xi_0(\cd)$. Note that due to the presence of the (strictly positive) impulse cost, the trivial impulse control is different from the zero impulse control (which contains impulses with $\xi_k=0$). All of these conditions are standard assumptions for impulse control problem.

\ms

We now collect some standard results into the following proposition.

\bp{state X} \sl Let {\rm(H1)$'$--(H1)$''$} hold. Then for any $(t,x)\in[0,T]\times\dbR^n$ and $\xi(\cd)\in\sK[t,T]$ of form \rf{xi}, state equation \rf{state1} admits a unique solution $X(\cd)=X(\cd\,;t,x,\xi(\cd))$, and the following estimates hold:
\bel{|X(s)|}|X(s)|\les e^{L(s-t)}(1+|x|)+|\xi(s)|+L\int_t^se^{L(s-\t)}|\xi(\t)|d\t,\qq s\in[t,T],\ee
%
%
If $\h X(\cd)=X(\cd\,;t,\hat x,\xi(\cd))$ with $\hat x\in\dbR^n$, then
\bel{|X-X|}|X(s)-\h X(s)|\les e^{L(s-t)}|x-\hat x|,\qq s\in[t,T].\ee
In addition, if {\rm (H2)} hold, then there exists a continuous increasing function $\h\n:[0,\i)\to[0,\i)$ such that the value function $V(\cd\,\cd)$ satisfies the following estimates:
\bel{V-V<L}\left\{\2n\ba{ll}
\ds0\les V(t,x)\les \h\n(|x|),\\
\ns\ds|V(t,x)-V(\,\h t,x)|\les\h\n\big(|x|\vee|\h x|\big)[|t-\h t\,|+|x-\h x|^\d]\qq\forall (t,x),(\h t,\h x)\in[0,T]\times\dbR^n.\ea\right.\ee
Moreover, Bellman's principle of optimality holds, i.e., for any $(t,x)\in[0,T)\times\dbR^n$, the following principle of optimality holds:
\bel{optimality-1}
V(t,x)\les\min_{\xi\in K}\big\{V(t,x+\xi)+\ell(t,\xi)\big\}\equiv\BN[V](t,x),\qq\forall(t,x)\in[0,T)\times \dbR^n,\ee
\bel{optimality-2}
V(t,x)\les\int_t^{\h t}g(r,X(r;t,x,\xi_0(\cd)))dr+V(\,\h t,X(\h t;t,x,\xi_0(\cd))),\q \forall0\les t\les\h t\les T,~x\in\dbR^n,\ee
and if the strict inequality holds in \rf{optimality-1}, then there exists
a $\bar t\in(t,T]$ such that
\be\label{t,bart,hatt}
V(t,x)=\int_t^{\h t}g(r,X(r;t,x,\xi_0(\cd)))dr+V(\,\h t,X(\h t;t,x,\xi_0(\cd))),\qq  0\les t\les\h t<\bar t\les T,~x\in\dbR^n.
\ee
Consequently, $V(\cd\,,\cd)$ is a viscosity solution of the HJB quasi-variational inequality
\bel{HJB-1}\left\{\ba{ll}
\ds\min\Big\{V_t(t,x)+\lan V_x(t,x),f(t,x)\ran+g(t,x),\BN[V](t,x)-V(t,x)\Big\}=0,\q
(t,x)\in[0,T)\times\dbR^n,\\
\ns\ds V(T,x)=h(x),\qq x\in\dbR^n.\ea\right.\ee
\ep

\ms

Finally, by the comparison theorem, we obtain that the value function $V(\cd,\cd)$ of Problem (C) is the unique viscosity solution to \rf{HJB-1}.

\section{Concluding Remarks}

We introduce a modification of the definition for viscosity super solution of the first order HJB type quasi-variational inequality. With such a modification, one has the comparison between viscosity super- and sub-solutions. It is not hard to see that for the second HJB type quasi-variational inequality, we can do the same thing, namely, modifying the definition of viscosity super-solution so that one can prove a comparison theorem. We prefer to omit the details here.

\baselineskip 18pt
\renewcommand{\baselinestretch}{1.2}

\end{document}